\documentclass[10pt,reqno]{article}
\usepackage{amsmath, amsthm, amssymb, stmaryrd}
\usepackage{hyperref}
\usepackage{enumerate}
\usepackage{url}
\usepackage{color}
\usepackage{tikz}
\usetikzlibrary{patterns}
\usepackage[T1,T2A]{fontenc}
\usepackage[utf8]{inputenc}

\theoremstyle{definition}

\theoremstyle{remark}

\def \le {\leqslant}
\def \ge {\geqslant}

\usepackage{tikz} 
\usetikzlibrary{calc,arrows,decorations.pathreplacing,fadings,3d,positioning}

\usepackage[mathscr]{euscript}

\begin{document}
 \begin{Large}
\centerline{On the number of  perfect positive definite 
quadratic forms}
\end{Large}

\vskip+0.3cm
\centerline{Nikolay Moshchevitin\footnote{Technische Universität Wien}}
\vskip+1cm

\begin{small}
{\bf Abstract:}
 We prove a new  asymptotic  lower bound for the number of perfect positive definite quadratic forms in $n$ variables which is close to the optimal one.
 \end{small}

\vskip+0.5cm
{\bf 1. Introduction.}

\vskip+0.3cm

Let 
$$
Q({\bf x}) =
\sum_{i,j=1}^n q_{i,j} x_ix_j,\,\,\,\,\, q_{i,j} = q_{j,i}
$$
be a positive definite quadratic form in $n$ real variables $ x_1,...,x_n$.
Two forms 
$Q'({\bf x})$ and $Q''({\bf x})$
are called {\it (integer) equivalent } if there exists a   $n\times  n$ matrix $\frak{U}$ with  integer entries and determinant
${\rm det}\, \frak{U} = \pm 1$ such that 
$Q''({\bf x})=Q'(\frak{U}{\bf x})$ for all ${\bf x}\in \mathbb{Z}^n$.

Let 
$$
m_Q = \min_{{\bf x} \in \mathbb{Z}^n\setminus\{{\bf 0}\}} Q({\bf x})
$$
be  arithmetic minimum  of the form $ Q({\bf x})$ and 
$$
\Sigma_Q = \{{\bf x}\in \mathbb{Z}^n: \,\,\, Q({\bf x}) = m_Q\}
$$
be the set of all integer vectors where this minimum  is attained.
Voronoï \cite{V} defined  form $Q({\bf x}) $ to be {\it perfect} if coefficients $q_{i,j} , 1\le i,j \le n$ are uniquely determined by the value of arithmetic minimum $m$ and the set $\Sigma$, that is  given $m$ and $\Sigma$ the form 
$Q({\bf x})$
with the given value of $ m_Q =m$ and with given set $\Sigma_Q =\Sigma$ is unique.
Voronoï realised
that this notion is important in accordance with
the problem of packing of balls in Euclidean space and proved that for any dimension $n$ the 
number of pairwise  non-equivalent perfect forms in $n$ variables is finite.
\vskip+0.3cm

We recall another related definition. Form $Q({\bf x})$ is called {\it extremal} if  any preserving determinant small perturbation of its coefficients $q_{i,j}$ 
makes the arithmetical minimum of the form  $Q({\bf x})$ smaller.

\vskip+0.3cm
It is well known that 
 \begin{equation}\label{qu}
\text{if}\,\,\,\, Q({\bf x}) \,\,\,\text{is extremal, then }
\,\,\, Q({\bf x})\,\,\, \text{is perfect}
\end{equation}
(see, for example \cite{Baran} Ch. II, \S 7.3).
Let $P_n$  and $E_n$ be the maximal number of pairwise non-equivalent positive definite quadratic forms in $n$ variables and extremal  quadratic forms in $n$ variables. From (\ref{qu}) it follows that  
$$E_n\le P_n.$$
The exact values of $P_n, E_n$ are known in small dimensions, namely
$$
P_1=P_2=P_3 = 1, \,\,\, P_4=2,\,\, \, P_5=3 ,\,\,\,  P_6 =7,\,\,\, P_7 = 33,\,\,\, P_8= 10916;
$$
$$
E_1=E_2=E_3 = 1, \,\,\, E_4=2,\,\, \, E_5=3 ,\,\,\,  E_6 =6,\,\,\, E_7 = 30,\,\,\, E_8= 2408,
$$
see discussion and exact references in \cite{VA}.

\vskip+0.3cm
As for the asymptotic bounds, to our knowledge the best known upper bound 
\begin{equation}\label{upe}
P_n < \exp\left( O(n^2\log n)\right)
\end{equation}
was proven by  vanWoerden \cite{W}
and  the lower bound 
 \begin{equation}\label{une}
P_n > \exp\left(n^{1-\varepsilon}\right)\,\,\,\,
\end{equation}
 for any positive $\varepsilon>0$  and for all $n$ large enough was obtained by
 Bacher \cite{B}.

\vskip+0.3cm
In the present paper we improve on the lower bound (\ref{une}) and show that the upper bound (\ref{upe}) is in a certain sense optimal.
Our result is as follows.

\vskip+0.3cm

{\bf Theorem 1.} {\it When $ n\to \infty$  one has asymptotic lower bound
$$
 E_n \ge  \exp\left(\frac{n^2\log n}{2} - O(n^2)\right).
$$}

 \vskip+0.3cm
 
 The rest of the paper contains a proof of Theorem 1.
 \vskip+0.3cm

{\bf 2. Lattice ${\bf A}_n$.}

\vskip+0.3cm
 Consider Euclidean space 
$\mathbb{R}^{n+1}   ({\bf z})$ with coordinates $ {\bf z} = (z_1,...,z_n,z_{n+1})$ and basis
$$
{\bf e}_1 = (1,0,...,0)^\top, ... , {\bf e}_{n+1} = (0,0,...,1)^\top.
$$
Geometric definition of $n$-dimensional lattice  ${\bf A}_n$ is as follows.  Consider
hyperplane
$$
\mathcal{H} = \{ {\bf z} = (z_1,...,z_n,z_{n+1})\top \in \mathbb{R}^{n+1}:\,\,\,\,\, z_1+...+z_n+z_{n+1} =0\} \subset\mathbb{R}^{n+1}
$$
  and define
$$
 {\bf A}_n  =\frac{1}{\sqrt{2}}\cdot  \mathbb{Z}^{n+1} \cap \mathcal{H}
$$
(it would be convenient for us to have in this definition factor $\frac{1}{\sqrt{2}}$  to make the length of the minimal vectors in ${\bf A}_n$ equal to 1).
Then $ {\bf A}_n$ is a complete lattice in $\mathcal{H}$ with a basis
\begin{equation}\label{gram}
{\bf f}_j =
\frac{1}{\sqrt{2}} ({\bf  e}_j - {\bf e}_{j+1}),\,\,\,\,\,\, j = 1, ...,n
\end{equation}
consisting of vectors of unit length
and covolume 
\begin{equation}\label{covol}
{\rm covol }\, {\bf A}_n = \frac{\sqrt{n+1}}{2^{n/2}}.
\end{equation}
So lattice ${\bf A}_n$ may be written as
$$
{\bf A}_n =\{ {\bf z} = x_1 {\bf f}_1+...+x_n{\bf f}_n,\,\,\,\, x_1,...,x_n \in \mathbb{Z}\}.
$$
Gram matrix of the collection (\ref{gram}) is of the form
$$
\frak{G} =
 \left(
 \begin{array}{cccccc}
 1&-\frac{1}{2}& \cdots &0&0&0\cr
  -\frac{1}{2}&1& \cdots &0&0&0\cr
   \vdots&\vdots& \vdots&\cdots  &\vdots&\vdots\cr
    0&0& \cdots&-\frac{1}{2} &1&-\frac{1}{2}\cr
     0&0& \cdots &0 &-\frac{1}{2}&1 \end{array}
 \right),\,\,\,\,\,
 {\rm det}\, \frak{G}
=\frac{n+1}{2^n}.
$$
It is clear that there exists a $n\times n$ matrix $\frak{A}$ such that 
$$
\frak{G} = \frak{A}^\top \frak{A}
 $$
 (the explicit expression for elements of $\frak{A}$ is not of importance for us).

 We identify hyperplane $ \mathcal{H} \subset \mathbb{R}^{n+1} $ with Euclidean space
 $\mathbb{R}^n   ({\bf x})$ supplied  with coordinates $ {\bf x} = (x_1,...,x_n)$ with respect to 
 an
orthogonal basis  of unit vectors 
\begin{equation}\label{basa}
{\bf g}_1 , ... , {\bf g}_{n} \in \mathbb{R}^n. 
\end{equation}
So lattice  $ {\bf A}_n$ may be identified with the lattice 
 \begin{equation}\label{base}
 \frak{A} \mathbb{Z}^n
 =
 \{{\bf x}= x_1 \frak{a}_1+...+x_n \frak{a}_n,\,\,\,\, x_j \in \mathbb{Z}\},
 \end{equation}
 where 
 $$
 \mathbb{Z}^n
 =\{{\bf x}= x_1 {\bf g}_1+...+x_n {\bf g}_n,\,\,\,\, x_j \in \mathbb{Z}\}
 $$ is the standard integer lattice with respect
 to basis (\ref{basa}) and
 $\frak{a}_j =\frak{A} {\bf g}_j$ are the columns of matrix $\frak{A}$, each of them has unit length:
 $$
 |\frak{a}_j|= 1,\,\,\,\,\, j = 1,..
 ..,n.
 $$
 Observe that 
 \begin{equation}\label{para}
 \text{diameter of fundamental parallelepiped}\,\,\, \Pi = 
 \{ {\bf x} = x_1 \frak{a}_1+...+x_n\frak{a}_n:\,\, \, 0\le x_j < 1,\,\,\, j=1,...,n\}\,\,
 \text{is less than}\,\,\, {n}.
 \end{equation}
 
  \vskip+0.3cm
 We recall some packing properties of lattice ${\bf A}_n$.
 
 \vskip+0.3cm
 \noindent
 {\rm({\bf i})}
 The interior of the unit ball 
 $$
 \Omega=\{{\bf x}\in \mathbb{R}^n:\,\, |{\bf x}| \le1\}
 $$
  contains no  non-zero points of ${\bf A}_n$.
  
   \vskip+0.3cm
 \noindent
 {\rm({\bf ii})}
  The only points of ${\bf A}_n$ which lie on the boundary of $\Omega$ are just  $n(n+1)$ points
   \begin{equation}\label{veee}
 \pm(\frak{a}_i+\frak{a}_{i+1} +...+ \frak{a}_j),\,\,\,\,\,\,\,\, 1\le i\le j \le n.
 \end{equation}
  
   \vskip+0.3cm
 \noindent
 {\rm({\bf iii})} Lattice ${\bf A}_n$ is {\it extremal}, that is  for any lattice $\Lambda\neq {\bf A}_n$ of the same covolume 
 ${\rm covol}\, \Lambda = {\rm covol}\,{\bf A}_n$   which is 
 close enough to ${\bf A}_n$ in the space of lattices modulo orthogonal transformations,   there exist  non-zero point ${\bf x}\in 
 \Lambda $  in the interior of $\Omega$.
  This property   is in fact  Theorem 7.1 from \cite{Baran}.
 
 \vskip+0.3cm
 \noindent
 {\rm({\bf iv})} For $n \ge 2$ 
lattice ${\bf A}_n$ has exactly $2(n+1)!$ different isometries, that is 
\begin{equation}\label{iso}
|\{ \frak{G}\in {\rm O}_n (\mathbb{R}) :\,\, \frak{G} {\bf A}_n = {\bf A}_n\}| = 2 (n+1)!
\end{equation}
(see \cite{ConveyS}, Ch 4, Section 6.1).

 \vskip+0.3cm
 
 The properties described above are proved and discussed in detail in \cite{Baran}   (see Ch. 2 \S 7) in terms of the corresponding 
  {\it Voronoi's first perfect positive definite quadratic form} $ \varphi ({\bf x})$. To help the reader we repeat few very basic properties of lattice ${\bf A}_n$
  in terms of form $ \varphi ({\bf x})$.
 One of the definitions of $ \varphi ({\bf x})$ is
 $$
 \varphi ({\bf x}) =
 \sum_{j=1}^n x_j^2 - \sum_{j=1}^{n-1} x_jx_{j+1} = {\bf x}^\top \frak{G} {\bf x} = |\frak{A}{\bf x}|^2. 
 $$
 Here we assume that coordinates  $x_1,...,x_n$ correspond to basis (\ref{basa}) in $\mathbb{R}^n$. In other words
 $$
  \varphi ({\bf x}) = \varphi_0 (\frak{A}{\bf x})\,\,\,\,\,\text{where}\,\,\,\,\,
    \varphi_0 ({\bf x}) =x_1^2+...+x_n^2.
 $$
 Properties  {\rm({\bf i})}, {\rm({\bf ii})} mean that  the 
 arithmetical minimum of the form $
 \varphi ({\bf x})$ is
 $$
 \min_{{\bf x}\in \mathbb{Z}^n \setminus\{{\bf 0}\}}  \varphi ({\bf x}) = 1,
 $$
 and it  attains at $n(n+1)$ vectors
 \begin{equation}\label{vee}
 \pm({\bf g}_i+ {\bf g}_{i+1}+...+  {\bf g}_j),\,\,\,\,\,\,1\le i \le j \le n 
 \end{equation}
 which  are the only integer vectors on the boundary of 
 ellipsoid
 $$
 \mathcal{E}=
 \{{\bf x}\in \mathbb{R}^n:\,\,\,\,\, \varphi ({\bf x}) \le 1\}.
 $$
 The set of these vectors  (\ref{vee}) we denote as $\Sigma_\varphi$.
Of course, the interior of  ellipsoid $
 \mathcal{E}$
 contains no  non-zero  points of integer lattice  $\mathbb{Z}^n$.
 
 Extremality of lattice ${\bf A}_n$  (property {\rm({\bf iii})}) is the same property as the  extremality of   quadratic form $\varphi ({\bf x})$.


  \vskip+0.3cm

 We need to say a few words about lattice  
 $${\bf A}_n^*
 =
 \{ {\bf x}\in\mathbb{R}^n:\,\,\, 
  ({\bf x}, {\bf y}) \in \mathbb{Z}\,\,\,\forall\, {\bf y} \in {\bf A}_n\}
 $$
  dual to ${\bf A}_n$.  
  It is obvious that  
  $ 2{\bf A}_n\subset {\bf A}_n^*$. 
  We should note that all the  entries of the inverse matrix $\frak{G}^{-1}$ are rational numbers with denominator $n+1$.
As 
$ (\frak{A}^{-1})^\top  = \frak{A}\frak{G}^{-1}$  
 we have 
 $$
 (n+1) {\bf A}_n^* = (n+1)  (\frak{A}^{-1})^\top \mathbb{Z}^n \subset\frak{A}\mathbb{Z}^n =  {\bf A}_n.
 $$
We conclude that 
\begin{equation}\label{1d0}
 (n+1) {\bf A}_n^* \subset {\bf A}_n\subset \frac{1}{2}  {\bf A}_n^*
 ,
 \end{equation}
and for indices we have
\begin{equation}\label{1d}
[{\bf A}_n:(n+1) {\bf A}_n^*]\,\, | \,  
\left[ \frac{1}{2}{\bf A}_n^*:(n+1) {\bf A}_n^*\right]\
=
(2(n+1))^{n}.
\end{equation}

Some other facts about the dual lattice ${\bf A}_n^*$  can be found  in \cite{ConveyS}, Ch. 4, Section 6.6.

 \vskip+0.3cm
{\bf 3. Sublattices and Abelian groups.}
 \vskip+0.3cm
Let  $p>n+1$ be prime and 
 \begin{equation}\label{aa}
 {\bf a}= (a_1,...,a_{n-1}, 1)^\top\in \mathbb{Z}^{n},\,\,\,\,\, 0\le a_j \le p-1, \,\,\ j=1,...,n-1
 \end{equation}
 be an integer vector.
 Define $\mathscr{W}$ to be the set of all $p^{n-1}$ integer vectors of the form (\ref{aa}).

  \vskip+0.3cm
 Define $n\times n$ matrix 
\begin{equation}\label{0}
   \frak{B}_{\bf a} = 
 \left(
 \begin{array}{ccccc}
 1&0& \cdots &0&\frac{a_1}{p}\cr
  0&1& \cdots &0&\frac{a_2}{p}\cr
   \vdots&\vdots& \vdots &\vdots&\vdots\cr
    0&0& \cdots &1&\frac{a_{n-1}}{p}\cr
     0&0& \cdots &0&\frac{1}{p}
 \end{array}
 \right),\,\,\,\,\,
  {\rm det}\,     \frak{B}_{\bf a}  = \frac{1}{p}
 \end{equation}
 and lattice
 $$
 \mathcal{L}_{\bf a} = \frak{A}_{\bf a}   \mathbb{Z}^n,\,\,\,\,\,
 \frak{A}_{\bf a}=
 \frak{A}   \frak{B}_{\bf a} 
 $$
 where matrix $\frak{A}$ is defined in the previous section.
 So $\mathcal{L}_{\bf a}$ has basis
 \begin{equation}\label{fraK}
 \frak{a}_1,...,\frak{a}_{n-1},  \frak{a}_n' =  \frac{ \frak{A} {\bf a}}{p} =\frac{a_1\frak{a}_1+...+a_{n-1}\frak{a}_{n-1} +\frak{a_n}}{p},
 \end{equation}
 where vectors $\frak{a}_j$ are defined in (\ref{base}),
 Then for any ${\bf a}$ lattice ${\bf A}_n$ will be a sublattice of lattice  $
 \mathcal{L}_{\bf a} $ of index
  \begin{equation}\label{mana0}
 [
 \mathcal{L}_{\bf a} :{\bf A}_n] = p
 .
 \end{equation}
  We should note that in particular
 \begin{equation}\label{aL}
 \frak{a}_n' =
 \frac{ \frak{A} {\bf a}}{p} \in \mathcal{L}_{\bf a}.
 \end{equation}
 Moreover
  \begin{equation}\label{section}
 \bigcup_{{\bf a} \in \mathcal{W} }\mathcal{L}_{\bf a} \subset \frac{1}{p} {\bf A}_n =\{ {\bf x} \in \mathbb{R}^n:\,\,\,\, p{\bf x} \in {\bf A}_n\}\,\,\,\,\,\,
 \text{and}
 \,\,\,\,\,\,\,
 \mathcal{L}_{\bf a} \cap  \mathcal{L}_{{\bf a}'}  = {\bf A}_n \,\,\, \forall \, {\bf a} \neq {\bf a}'.
 \end{equation}

 Let 
 $$
 \mathcal{L}_{\bf a}^* = \{ {\bf x}\in \mathbb{R}^n: \,\,\, ({\bf x}, {\bf y}) \in \mathbb{Z}\,\,\,\forall\, {\bf y} \in \mathcal{L}_{\bf a}\}
 $$ be the lattice dual to $ \mathcal{L}$.
 Notice that
 by the definition and (\ref{aL}) we see that 
  \begin{equation}\label{aL1}
   {\bf x} \in  \mathcal{L}_{\bf a}^* \,\,\,\,\,\,\,\,\Longrightarrow\,\,\,\,\,\,\,\,
   \left({\bf x},  \frak{a}_n'  \right) \in \mathbb{Z}.
   \end{equation}
   Moreover 
   $$
   \mathcal{L}_{\bf a}^* \subset {\bf A}_n^*\,\,\,\,\,
   \text{and}
   \,\,\,\,\,
   [{\bf A}_n^*:   \mathcal{L}_{\bf a}^*]=
    [
 \mathcal{L}_{\bf a}:{\bf A}_n] = p.
   $$
     Taking into account  (\ref{1d0}, \ref{1d})   we conclude that 
  \begin{equation}\label{mana3}
 (n+1)  \mathcal{L}_{\bf a}^* \subset   (n+1){\bf A}_n^* \subset {\bf A}_n
     \end{equation}
 and
 \begin{equation}\label{mana1}
 p =
  [{\bf A}_n^*:\mathcal{L}_{\bf a}^* ]=
   [(n+1){\bf A}_n^*: (n+1)\mathcal{L}_{\bf a}^* ]
 \, |
 \,
  [{\bf A}_n:  (n+1)  \mathcal{L}_{\bf a}^*] 
  =
    [{\bf A}_n:  (n+1)  {\bf A}_n^*] 
  \cdot 
    [(n+1){\bf A}_n^*:  (n+1)  \mathcal{L}_{\bf a}^*] 
  \, |\, p(2(n+1))^n.
    \end{equation}
  So finally
  \begin{equation}\label{mana}
   (n+1)  \mathcal{L}_{\bf a}^* \subset {\bf A}_n \subset \mathcal{L}_{\bf a},\,\,\,\,\,
  \end{equation}
  and the corresponding indices satisfy (\ref{mana0}, \ref{mana1}). We see that  the factor $\mathscr{X} = \mathcal{L}_{\bf a}/(n+1)  \mathcal{L}_{\bf a}^* $
  is an Abelian group of cardinality
  $$
  |\mathscr{X} | = [\mathcal{L}_{\bf a}:(n+1)  \mathcal{L}_{\bf a}^*]=
 [\mathcal{L}_{\bf a}: {\bf A}_n] \cdot [{\bf A}_n:(n+1)  \mathcal{L}_{\bf a}^*]
  $$
  which satisfies 
  $$
 p^2\, | \,  |\mathscr{X} |\, |\, p^2 (2(n+1))^n,
  $$
  and  another  factor group $ \mathscr{Y} = \mathcal{L}_{\bf a}/{\bf A
  _n} $ of cardinality 
  $$
  |\mathscr{Y} | =   [
 \mathcal{L}_{\bf a} :{\bf A}_n] = p
.
  $$
   \vskip+0.3cm
 {\bf Lemma 1.}   {\it  Let $p\, \not | \, 2(n+1)$ and
 \begin{equation}\label{phi}
 \varphi ({\bf a}) \not\equiv 0 \pmod{p}.
 \end{equation}
 Then there is the unique subgroup of $\mathscr{X}$ of cardinality $p$ as well as the unique subgroup of index $p$.}
\vskip+0.3cm
 
 {\bf Corollary.}   {\it  If ${\bf a}$ satisfies (\ref{phi}) and sublattice
 $A$ satisfies 
 \begin{equation}\label{cora}
   (n+1)  \mathcal{L}_{\bf a}^*\subset  A \subset \mathcal{L}_{\bf a}
 \,\,\,\,\,\text{and}
 \,\,\,\,\,
  [\mathcal{L}_{\bf a}: A] = p,
  \end{equation}
   then
 $A= {\bf A}_n$.}
\vskip+0.3cm

Proof of Lemma 1.
First of all we show that there exists a cyclic subgroup $\mathscr{X}_1 \subset \mathscr{X}$ of cardinality $p^2$. Indeed, if it is not the case,
as
 $(p, 2(n+1)) = 1$,
Abelian group $\mathscr{X}$  can be written as
$$
\mathscr{X} = \mathscr{X}_1' \oplus \mathscr{X}_1''\oplus \mathscr{X}_2\,\,\,\,\,
\text{where}\,\,\,\,\,
|\mathscr{X}_1'| =|\mathscr{X}_1''|  = p,\,\,\,\,\,
|\mathscr{X}_2|  \, |\,  (2(n+1))^n.
$$
Then for any $\overline{\frak{a}}\in \mathscr{X}$ the element $p(2(n+1))^n \overline{\frak{a}}$ will be zero element in $\mathscr{X}$, and  so
 \begin{equation}\label{phi0}
p(2(n+1))^n {\frak{a}}\in  (n+1) \mathcal{L}_{\bf a}^*\subset \mathcal{L}_{\bf a}^*.
 \end{equation}

In  $\mathscr{X}$ consider element 
$$
\overline{\frak{a}_n'} =  \frak{a}_n'+ (n+1)  \mathcal{L}_{\bf a}^* \in \mathscr{X},
$$
where $\frak{a}_n'$ is defined in (\ref{aL}).
We see that 
$$
p(2(n+1))^n \frak{a}_n' =  \frak{A}(2(n+1))^n {\bf a} +    (n+1) \mathcal{L}_{\bf a}^*
,
$$
{where}
$
(2(n+1))^n{\bf a} $ 
{is an integer vector}.
But  (\ref{phi0})  allows to apply (\ref{aL1}) to ${\bf x} = \frak{A}(2(n+1))^n {\bf a} \in \mathcal{L}_{\bf a}^*$. In such  a way we get
 $$
 (\frak{A}(2(n+1))^n {\bf a},\frak{a}_n') =
 (2(n+1))^n  (\frak{A} {\bf a}, \frac{\frak{A} {\bf a}}{p}) \in \mathbb{Z},
 $$
 and so
 $$
 \varphi ({\bf a})  =   (\frak{A} {\bf a},\frak{A} {\bf a})\equiv 0\pmod{p},
$$
and this is a contradiction to condition (\ref{phi})  of Lemma 1.  

\vskip+0.3cm
We proved that $\mathscr{X}$ contains a cyclic subgroup of order $p^2$. So
$$
\mathscr{X} = \mathscr{X}_1 \oplus \mathscr{X}_2,\,\,\,\,\,
\text{where}\,\,\, \mathscr{X}_1\,\,\,\text{is a cyclic group of order}\,\,\, p^2\,\,\,\text{ and}\,\,\, 
(|\mathscr{X}_2|, p)=1.
$$
Let $\mathscr{Y}' \subset \mathscr{X}$ be a subgroup of cardinality $p$.
Then
$$
\mathscr{Y}'  =
\mathscr{Y}' _1\oplus \mathscr{Y}'_2,\,\,\,\, \mathscr{Y}' _j \subset \mathscr{X}_j,\,\,\,\,\, 
|\mathscr{Y}' _1|\cdot | \mathscr{Y}'_2|=p.
$$
It is clear that $|\mathscr{Y}' _1|=p,  | \mathscr{Y}'_2|=1$. Moreover the subgroup of the order $p$ in $\mathscr{X} _1$ is unique.
The similar argument shows the uniqueness of subgroup  $\mathscr{Y}'' \subset \mathscr{X}$  of index $p$.

Lemma is proven.$\Box$
 
 \vskip+0.3cm
{\bf 4. Quadratic forms.}
 \vskip+0.3cm
 For any ${\bf a}$ we consider quadratic form
 $$
 \mathcal{Q}_{\bf a}({\bf x}) =\varphi (\frak{B}_{\bf a}{\bf x})  = {\bf x}^\top \frak{Q}_{\bf a} {\bf x},\,\,\,\,\,\,   \frak{Q}_{\bf a} =  (\frak{A}_{\bf a} )^\top \frak{A}_{\bf a} 
 .$$
 
   \vskip+0.3cm
 {\bf Lemma 2.}  
  {\it
  For any vector ${\bf a}\in \mathscr{W}$ 
  satisfying (\ref{phi})
   there exists not more than $ 2(n+1)!$  different vectors ${\bf a}'$ of the  form (\ref{aa}) such that the forms
  $ \mathcal{Q}_{\bf a}({\bf x})  $ and $ \mathcal{Q}_{{\bf a}'}({\bf x})$ are integer equivalent.
  }
     \vskip+0.3cm
     Proof. 
     Equivalence of the forms means that 
     $$
     \frak{U}^\top\frak{Q}_{\bf a} \frak{U} = \frak{Q}_{{\bf a} '}
     \,\,\,\,\,
     \text{
     for some }\,\,\,\,\, \frak{U} \in {\rm GL}_n (\mathbb{Z}).
     $$
     Consider matrix
     $$\frak{O} = \frak{A}_{\bf a} \frak{U} \frak{A}_{{\bf a}'}^{-1}.$$
     Then
     $$
     \frak{O}^\top\frak{O} = 
     (\frak{A}_{{\bf a}'}^\top)^{-1} \frak{U}^\top \frak{A}_{{\bf a}}^{\top}
     \frak{A}_{\bf a} \frak{U} \frak{A}_{{\bf a}'}^{-1} =  (\frak{A}_{{\bf a}'}^\top)^{-1}
     \frak{U}^\top\frak{Q}_{\bf a}
     \frak{U} \frak{A}_{{\bf a}'}^{-1}  =
      (\frak{A}_{{\bf a}'}^\top)^{-1} \frak{U}^\top \frak{A}_{{\bf a}}^{\top}
     \frak{A}_{\bf a} \frak{U} \frak{A}_{{\bf a}'}^{-1} =  (\frak{A}_{{\bf a}'}^\top)^{-1}
     \frak{Q}_{{\bf a}'}
     \frak{A}_{{\bf a}'}^{-1}
     = \frak{I}
     $$ 
     is the unit $n\times n$ matrix.
So $ \frak{O} \in {\rm O}_n (\mathbb{R}) $ and  $\frak{O}  \frak{A}_{{\bf a}'} = \frak{A}_{\bf a} \frak{U}  $, that is
$
\frak{O}\mathcal{L}_{{\bf a}'} = \mathcal{L}_{{\bf a}}$. Now we see that 
$$
A=  \frak{O} {\bf A}_n \subset  \mathcal{L}_{{\bf a}}
$$
 is a subgroup of index $p$. Moreover, from (\ref{mana}) we deduce that 
$$
(n+1) \mathcal{L}_{{\bf a}}^* =
(n+1)(\frak{O} \mathcal{L}_{{\bf a}'})^* =
\frak{O}  (n+1)  \mathcal{L}_{{\bf a}'}^* \subset \frak{O} {\bf A}_n = A,
$$
so $A$ satisfies (\ref{cora}).
By Corollary to Lemma 1 we conclude that $ \frak{O} {\bf A}_n = {\bf  A}_n$. Formula (\ref{iso}) finalises the proof.$\Box$.

    \vskip+0.3cm
 {\bf Lemma 3.}  
  {\it Assume that lattice $ \mathcal{L}_{\bf a}$ has no  points in the  closed unit ball $\Omega$ different from points of ${\bf A}_n$.
  Then  quadratic form $\mathcal{Q}_{\bf a}$ is extremal.}
  \vskip+0.3cm
  Proof. This lemma is almost obvious, because  vectors  (\ref{veee}) belong $ \mathcal{L}_{\bf a}$  and by the construction there is positive $\varepsilon$ such that 
  $$
  (1+\varepsilon) \Omega \cap  \mathcal{L}_{\bf a} =  (1+\varepsilon) \Omega \cap {\bf A}_n.
  $$
   This leads to 
  $$
  m_{\mathcal{Q}_{\bf a}} = m_{\varphi} = 1.
  $$
Let quadratic form  $\mathcal{Q}'({\bf x}) = {\bf x}^\top \frak{Q}' {\bf x} ,\,\, \frak{Q}' \neq \frak{Q}_{\bf a}$  be a 
  sufficiently small, nontrivial perturbation of 
  $\mathcal{Q}_{\bf a}({\bf x}) =
  \varphi (\frak{B}_{\bf a}{\bf x}) =
   {\bf x}^\top \frak{Q}_{\bf a} {\bf x} $ with the same determinant 
  ${\rm det}\, \frak{Q}' = {\rm det}\, \frak{Q}_{\bf a}$.
  Then form
  $\mathcal{Q}'( \frak{B}_{\bf a}^{-1}{\bf x}) \neq   \varphi ({\bf x})$
is a small nontrivial perturbation of   $\varphi ({\bf x})$ with the fixed determinant.
As $\varphi$ is extremal, we conclude that there exists ${\bf x}_0 \in \mathbb{Z}^n\setminus \{{\bf 0}\}$ such that 
$$
\mathcal{Q}'( \frak{B}_{\bf a}^{-1}{\bf x}_0) <1.
$$
But as 
 $$
    \frak{B}_{\bf a}^{-1} =  
 \left(
 \begin{array}{ccccc}
 1&0& \cdots &0& -{a_1}\cr
  0&1& \cdots &0&-{a_2}\cr
   \vdots&\vdots& \vdots &\vdots&\vdots\cr
    0&0& \cdots &1&-{a_n}\cr
     0&0& \cdots &0&{p}
 \end{array}
 \right) 
 $$
 is an integer matrix, we see that $ {\bf y}_0 =\frak{B}_{\bf a}^{-1}{\bf x}_0\in \mathbb{Z}^n$.
 So 
$$
m_{\mathcal{Q}'} \le \mathcal{Q}'( {\bf y}_0) <1,
$$
and this proves extremality of $\mathcal{Q}_{\bf a}$.$\Box$ 

   \vskip+0.3cm
    From (\ref{qu}) we immediately get 
    \vskip+0.3cm
 {\bf Corollary.}   {\it $\mathcal{Q}_{\bf a}$ is perfect.}
      \vskip+0.3cm
   
 {\bf Remark.}    We can easily observe perfectness of $\mathcal{Q}_{\bf a}$ directly.
 Indeed,
      $\varphi ({\bf x })$  is perfect. Consider   the set
     $
     \Sigma_\varphi \subset \mathbb{Z}^n
     $
      of all representations  of arithmetic minima 
      of $\varphi ({\bf x })$. Then the system of equations
     \begin{equation}\label{sys}
     {\bf w}^\top Q{\bf w} =
     \sum_{1\le i,j \le n } q_{i,j} w_i w_j =0\,\,\,\,\, q_{i,j} = q_{j,i};\,\,\,\,\,{\bf w} = (w_1,...,w_n) \in  \Sigma_\varphi
     \end{equation}
     for entries of matrix $ Q =(q_{i,j})_{i,j+1}^n$
     has only trivial solution $q_{i,j} = 0 \,\, \forall i,j$. 
     Consider $\Sigma ' = \frak{B}_{\bf a}^{-1} \Sigma_\varphi \subset \mathbb{Z}^n$.
  From the conditions of Lemma 3 we see that every $ {\bf w}'\in \Sigma '$ represents the arithmetical minimum of the form  $\mathcal{Q}_{\bf a}$ 
  which coincides with the arithmetical minimum of  $\varphi$.
 Instead of the  system of equations (\ref{sys}) we consider the system
 $$
  ({\bf w}')^\top Q' {\bf w}' = 0,\,\,\,\, (Q')^\top = Q';\,\,\,\,\,\,{\bf w}'\in  \Sigma'
  $$
  with unknown entries of real $n\times n$ matrix $Q'$. Write $Q=(\frak{B}_{\bf a}^\top)^{-1} Q'\frak{B}_{\bf a}^{-1}$.
  Then $Q$ satisfies (\ref{sys}) and so is zero matrix. But then $Q'$ is also a zero matrix. This means that  quadratic form 
   $\mathcal{Q}_{\bf a}$ is perfect.

    \vskip+0.3cm
{\bf 5. Bounds for the number of lattice points.}
 \vskip+0.3cm
 
 {\bf Lemma 4.}  
  {\it  The number  $N_p$ of pointы of the lattice $ \frac{1}{p}{\bf A}_n$ in the unit ball $\Omega$ is 
  $\le  \frac{2^{n/2}(p+{n})^n {\rm vol}\, \Omega}{\sqrt{n+1}}$.}
  
       \vskip+0.3cm
     Proof. If
  $$
  {\bf x} \in  \frac{1}{p}{\bf A}_p\cap\Omega , 
  $$
  by (\ref{para}) we see that 
  $$
    {\bf x}+ \frac{1}{p } \Pi \subset \left(1+\frac{{n}}{p}\right) \Omega.
    $$
    So
    $$
    \bigcup_{\bf{x} \in \frac{1}{p}{\bf A}_p}    \left({\bf x}+ \frac{1}{p }  \Pi  \right) \subset \left(1+\frac{{n}}{p}\right) \Omega.
    $$
    If ${\bf x}, {\bf x}' \in   \frac{1}{p}{\bf A}_p $ are different, we have
    $$
     \left({\bf x}+ \frac{1}{p }  \Pi  \right)\cap
     \left({\bf x}'+ \frac{1}{p }  \Pi  \right) = \varnothing.
     $$
    By considering volume we conclude that 
    $$
    \frac{N_p {\rm vol}\, \Pi }{p^n} \le \left(1+\frac{{n}}{p}\right)^n 
    {\rm vol}\, \Omega.
    $$
    Taking into account (\ref{covol}) proves lemma.$\Box$

     \vskip+0.3cm
     
     The following lemma uses standard mean value argument similar to that from Minkowski-Hlawka's theorem (\cite{Cas}, Ch. VI, see also  argument from \cite{Corobo}).
 
     \vskip+0.3cm
 {\bf Lemma 5.}  
  {\it  
  Assume that 
  \begin{equation}\label{conti}
  \frac{2^{n/2}\left(1+\frac{{n}}{p}\right)^n {\rm vol}\, \Omega}{\sqrt{n+1}}< \frac{1}{2p}.
  \end{equation}
   Then there exists a set
     $\mathscr{W}_1 \subset \mathscr{W} $   such that  $ |\mathscr{W}_1|=
\left[ p^{n-1}/2\right]$ and 
  $$\mathcal{L}_{\bf a}\cap  \Omega = {\bf A}_n\cap \Omega,\,\,\,\,\,
  \forall \, {\bf a} \in \mathscr{W}_1.
  $$
  }
  
       \vskip+0.3cm
     Proof.  
     We construct   set $\mathscr{W}_1$ inductively.
     Let  
     $$ 
      \chi ({\bf x}) 
     =\begin{cases}
     1,\,\,\,\,\,{\bf x} \in  \Omega,\cr
          0,\,\,\,\,\,{\bf x} \not\in  \Omega
          \end{cases}
          \,\,\,\,\text{and}
          \,\,\,\,\,\,\,\,
     \chi '({\bf x}) 
     =\begin{cases}
     1,\,\,\,\,\,{\bf x} \in  \Omega \setminus {\bf A}_n,\cr
          0,\,\,\,\,\,{\bf x} \not\in  \Omega \setminus {\bf A}_n
          \end{cases}
          \,\,\,\,\text{so}\,\,\,\,
           \chi' ({\bf x}) \le  \chi ({\bf x}) .
     $$ 
   Consider sum
     $$
     S({\bf a}) = 
     \sum_{{\bf x}\in \mathcal{L}_{\bf a} }\chi' ({\bf x}) = |\Omega\cap (\mathcal{L}_{\bf a} \setminus{\bf A}_n)|.
     $$
     Recall that  ${\bf A}_n$ does not have non-zero points in the interior of $\Omega$,
     while on the boundary of $\Omega$ it has $ n(n+1)$ points.
      Taking into account (\ref{section}) we continue with 
       $$
       \sum_{{\bf a} \in \mathscr{W}}
     S({\bf a}) \le \sum_{{\bf x}\in \frac{1}{p}{\bf A}_n} \chi' ({\bf x})  
     <
     \sum_{{\bf x}\in \frac{1}{p}{\bf A}_n} \chi ({\bf x}) 
     = N_p.
     $$
     Taking into account condition (\ref{conti})   and Lemma 4 we conclude that 
     $$
     \frac{1}{|\mathscr{W}|} \sum_{{\bf a} \in \mathscr{W}}      S({\bf a}) <
     \frac{N_p}{p^{n-1}} <\frac{1}{2}<1.
     $$
     This means that  there exists $ {\bf a}_1 \in \mathscr{W}$ such that 
     $     S({\bf a}_1) = 0$, that is  $\mathcal{L}_{{\bf a}_1} \cap \Omega = {\bf A}_n\cap \Omega.$
     
     Now we assume that we have constructed vectors ${\bf a}_1,..., {\bf a}_k\in \mathscr{W},
     \, k < |\mathscr{W}|/2-1$ such that  $\mathcal{L}_{{\bf a}_j} \cap \Omega = {\bf A}_n\cap \Omega$ for every $ j \le k$.
     Consider sum
            $$
       \sum_{{\bf a} \in \mathscr{W} \setminus \{ {\bf a}_1,..., {\bf a}_k\} }
     S({\bf a}) \le
            \sum_{{\bf a} \in \mathscr{W} }
     S({\bf a}) 
    < N_p.
     $$
     We see that  as $ k< |\mathscr{W}|/2 = p^{n-1}/2$, again
      $$
     \frac{1}{|\mathscr{W} \setminus \{ {\bf a}_1,..., {\bf a}_k\} |} \sum_{{\bf a} \in \mathscr{W} \setminus \{ {\bf a}_1,..., {\bf a}_k\} }      S({\bf a})
     <
     \frac{N_p}{p^{n-1}/2} <1,
     $$
     and
     this means that  there exists $ {\bf a}_{k+1} \in \mathscr{W} \setminus \{ {\bf a}_1,..., {\bf a}_k\} $ such that 
     $     S({\bf a}_{k+1}) = 0$, that is  $\mathcal{L}_{{\bf a}_{k+1}} \cap \Omega = {\bf A}_n\cap \Omega.$
     In such a way we construct the desired set $\mathcal{W}_1=
     \{ {\bf a}_1,..., {\bf a}_{\left[ p^{n-1}/2\right]}\}$.$\Box$

         \vskip+0.3cm
{\bf 6. Proof of Theorem 1.}
 \vskip+0.3cm
Consider set 
$$
\mathscr{W}_2 = \{ {\bf a}\in \mathscr{W}:\,\,\,\,\, \varphi ({\bf a}) \equiv 0\pmod{p}\}.
$$
Notice that for fixed $a_1,...,a_{n-2}$ there may be not more than two values of $a_{n-1}$ for which 
$ \varphi (a_1,....,a_{n-2}, a_{n-1},1) \equiv 0\pmod{p}$, and so 
\begin{equation}\label{ww}
|\mathscr{W}_2|\le 2p^{n-2}.
\end{equation}
For large $n$ 
take prime $p$ in the range 
\begin{equation}\label{ww1}
\frac{1}{2} \left( \frac{n}{2^5\pi e}\right)^{n/2} \le p 
\le 
 \left( \frac{n}{2^5\pi e}\right)^{n/2}.
\end{equation}
As 
\begin{equation}\label{ww2}
{\rm vol}\,\Omega = \frac{\pi^{n/2}}{\Gamma\left(\frac{n}{2}+1\right) } \sim
    \frac{1}{\sqrt{\pi n}}\left(\frac{2\pi e}{n}\right)^{n/2},\,\,\,\,\, n\to \infty,
     \end{equation}
     we see that (\ref{conti}) is  satisfied.
   
   Define 
    $
    \mathscr{W}_3 = \mathscr{W}_1 \setminus \mathscr{W}_2,
    $
    where $\mathscr{W}_1 $ comes from
    Lemma 5. By (\ref{ww}, \ref{ww1}, \ref{ww2}) we deduce that 
    $$
    | \mathscr{W}_3 | \ge  | \mathscr{W}_1 |  -  | \mathscr{W}_2 |   = \frac{p^{n-1}}{2}(1-o(1)) = \exp\left(\frac{n^2\log n}{2} - O(n^2)\right).
    $$
    By Lemma 3 we see that  for ${\bf a} \in \mathscr{W}_3$ the form $ \mathcal{Q}_{\bf a}$ is extremal.
    By Lemma  2 the number of ${\bf a}'$ with quadratic   form $ \mathcal{Q}_{{\bf a}'}$ equivalent to 
    $ \mathcal{Q}_{\bf a}$ is $O((n+1)!) = O (\exp (n\log n))$. So we have constructed
    $\exp\left(\frac{n^2\log n}{2} -  O(n^2) -O(n\log n)\right)$ non-equivalent extremal forms in $n$ variables.

   Theorem is proven.$\Box$
   
   \vskip+0.3cm
   {\bf Acknowledgements.}
   The author thanks Vasiliy Neckrasov for many fruitful discussions of the manuscript. This research is supported by 
  by Austrian Science Fund (FWF), Forschungsprojekt PAT1961524.


\begin{thebibliography}{1} 
 

 
\bibitem{B}
 R. Bacher,\,  On the number of perfect lattices, J. Théor. Nr. Bordeaux, 30:3 (2018), 917–945.

\bibitem{Cas}
J.W.S. Cassels,\,
An Introduction
to the Geometry
of Numbers, Springer, 1971.

\bibitem{ConveyS}
J. H. Conway, N. J. A. Sloane,\, 
Sphere Packings, Lattices and Groups,
Springer,
1999.

\bibitem{Corobo}
 N. M. Korobov,\,  Properties and calculation of optimal coefficients, Dokl. Akad. Nauk SSSR, 132:5 (1960), 1009–1012.
 
 
 \bibitem{Baran}
S. S. Ryshkov, E. P. Baranovskii, \, Classical methods in the theory of lattice packings, Russian Math. Surveys, 34:4 (1979), 1–68.
 
 
 \bibitem{VA}
 M.D. Sikirić, A. Schürmann, F. Vallentin, 
 Classification of eight-dimensional perfect forms,
 Electronic research announcements of AMS, 13 (2007), 21 - 32.
 

 
 \bibitem{V}
 G. Voronoï, \, Nouvelles applications des paramètres continus a la theorie des formes quadratiques. Deuxieme memoire. Recherches sur les parallelloèdres primitifs, J. Reine Angew. Math. 134 (1908), 198–287.

\bibitem{W}
W.P.J. van Woerden,\, 
An upper bound on the number of perfect quadratic forms,
Advances in Mathematics 365 (2020),  107031, 12p.
\end{thebibliography}
 \end{document}